\documentclass[12pt,a4paper]{amsart}
\usepackage{amsmath, amssymb, bm}
\usepackage[all]{xy}

\let \lan \langle
\let \ran \rangle

\def \c {\colon}

\def \q {\qquad}
\def \qq {\qquad \qquad}

\def\dd {\partial}
\def \sep {\,  |  \,}     

\def \skp {\medskip}
\def \ndt {\noindent}
\def \Ndt {\medskip  \noindent}

\def \tilde {\raise.17ex\hbox{$\scriptstyle\mathtt{\sim}$}}   

\def \ard {\ar@{-->}}
\def \arp {\ar@{.>}}
\def \arv {\ar@{}}   
\def \arl {\ar@{-}}    
\def \arld {\ar@{--}}    
\def \arlp {\ar@{..}}    

\def \sup {{\rm sup}}
\def \inf {{\rm inf}}
\def \ti {\! \times \!}
\def \te {\otimes}

\def \sq {{\, {\scriptstyle\square}\, }}

\def \Sum {\raisebox{0.45ex}{${\mbox{\fontsize{10}{10}\selectfont\ensuremath{\sum}}}$}}

\def \setm {{\raise.4ex\hbox{$ \, \scriptscriptstyle{\setminus} \; $}}}
\def \sub {\subset}
\def\le{\leqslant}
\def\ge{\geqslant}

\def \and {\mbox{ and }}
\def \for {\mbox{for }}

\def \inn {\mbox{ in }}
\def \IF {\; \mbox{ if } \:}

\def \al {\alpha}

\def \ga {\gamma}
\def \de {\delta}
\def \ep {\varepsilon}

\def \la {\lambda}
\def \si {\sigma}

\def \La {\Lambda}

\def \id {{\rm id\,}}
\def \Ob {{\rm Ob\,}}

\def \Ut {{\underline{t}}}

\def \Uga {{\underline{\gamma}}}
\def \Vv {{\vec{v}}}
\def \Vw {{\vec{w}}}
\def \Bw {{\bf w}}
\def \Bg {{\bf g}}
\def \Bq {{\bf q}}
\def \Oinf {\overline{\infty}}

\def \B {\mathsf{B}}
\def \C {\mathsf{C}}

\def \Set {\mathsf{Set}}
\def \Ab {\mathsf{Ab}}
\def \Ord {\mathsf{Ord}}
\def \Mtr {\mathsf{Mtr}}

\def \Cat {\mathsf{Cat}}
\def \one {{\bf 1}}
\def \two {{\bf 2}}
\def \wPi {{\rm w} \Pi}
\def \bbX {\mathbb{X}}

\def \es {\varnothing}
\def \sing  {\{*\}}
\def \bbZ {\mathbb{Z}}
\def \bbI {\mathbb{I}}
\def \bbR {\mathbb{R}}
\def \bR {{\bf R}}
\def \OR {\overline{\mathbb{R}}}

\begin{document}

\title[Real metrics and Lorentz manifolds]{Enriched categories, real metrics and Lorentz manifolds}

\author[M. Grandis]{Marco Grandis}

\address{Marco Grandis, Dipartimento di Matematica, Universit\`a di Genova, 16146-Genova, Italy.}
\email{grandismrc@gmail.com}

\subjclass{18D20, 83Axx, 54E35}

\keywords{Enriched categories, Lorentz manifold, Minkowski spacetime, generalised metric space}

\begin{abstract}

This expository article brings together two subjects: generalised metrics based on enriched 
categories, on the one hand, and Lorentz manifolds, on the other, at the price of dealing with 
details that are well known {\em either} in category theory {\em or} in relativity.

	The spacetime of relativity can be given a {\em real valued metric} $\rho(x, y)$, with values 
in the extended real line $\OR= [- \infty, \infty]$, or better (if equivalently) a real valued 
`antimetric' $\ga(x, y) = - \rho(x, y)$ (satisfying a reverse triangle inequality); the latter, as a 
function of $y$, is positive on the timecone of $x$, annihilates on its lightcone, and is $- \infty$ 
on all events which cannot be influenced by $x$.

	All this can be given a well-established base in category theory, extending Lawvere's 
notion of a metric space. In fact, a space with a real valued metric can be viewed as an 
enriched category on the extended real line $ \OR$, structured as a symmetric monoidal 
closed category.

\end{abstract}

 \maketitle

\section*{Introduction}\label{Intro}

\subsection{Real valued metrics}\label{0.1}
Lawvere's article \cite{L1} presented, in 1974, a (generalised) metric space $ X $ as a small 
category enriched on the positive extended half-line $ \Bw^+ = \OR^+ = [0, \infty]$; the latter is 
viewed as a category (with morphisms $ \la \ge \mu$) endowed with a strict symmetric 
monoidal closed structure given by the (extended) sum $ \la + \mu$.

	Concretely, $ X $ is a set equipped with a generalised metric $ \de\c X \ti X$ $\to \Bw^+ $ 
(the enriched hom) satisfying the triangular inequality (enriched composition) and a second 
condition (enriched units)
    \begin{equation}
\de(x, y) + \de(y, z)   \ge  \de(x, z),   \q\;\;   0   \ge  \de(x, x)   \qq   (x, y, z \in X).
    \label{0.1.1} \end{equation}

	The latter amounts to $ \de(x, x) = 0$, for all $ x \in X$. Symmetry and `separation' are 
not required.

	Although various similar extensions have been studied in the theory of metric spaces, 
Lawvere's structure combines a formal elegance (or even `necessity') within category theory, 
with concrete advantages, e.g.\ with respect to the existence of infinite products and 
coproducts -- which obviously fails for metrics with finite values.

	All this was recently extended in \cite{G4}, studying `real valued metrics', to express 
`profits' and `losses', in any sense -- possibly related to a variable of any science, from Physics 
to Economy. The positive extended half-line $ \Bw^+ $ is replaced by the whole extended real 
line $ \Bw = \OR= [- \infty, \infty]$. This complete lattice is again viewed as a category, with 
morphisms $ \la  \ge \mu$, and endowed with a strict symmetric monoidal closed structure 
given by the extended sum $ \la + \mu$. (Here we let $ - \infty + \infty = \infty$, because the
{\em initial} object $\infty$ must be preserved by adding any element.)

	A $\rho$-{\em metric space}, or $\rho$-{\em space}, is thus a small category $ X $ enriched 
on $ \Bw$. This means a set $ X $ with a mapping $ \rho\c X \ti X \to \OR$, called a 
{\em real metric}, or $\rho$-{\em metric}, that satisfies the same conditions as $ \de$, above
    \begin{equation}
\rho(x, y) + \rho(y, z)   \ge  \rho(x, z),   \q\;\;   0   \ge  \rho(x, x)   \qq   (x, y, z \in X),
    \label{0.1.2} \end{equation}
but the second means now that $\rho(x, x)$ can only be 0 or $ - \infty$.

	The extended real number $ \rho(x, y) $ expresses the cost of the transition from $ x $ to 
$ y$. Equivalently, the `antimetric' $ \ga(x, y) = - \rho(x, y) $ represents its gain; of course, it 
satisfies a reverse triangular inequality (as it appears in Relativity). The transition is 
unaffordable, or forbidden, when its cost is $ \infty$, and the gain is $ - \infty$.

\subsection{Outline and notation}\label{0.2}
Section 1 reviews the general notion of enriched category, dealing mostly with enrichment over 
a {\em strict} monoidal category. Section 2 shows how generalised metrics with positive or 
relative real values can be formalised in this way. The last two sections show how the Minkowski 
spacetime, or more generally a Lorentz manifold, can be given a real valued antimetric which 
forbids a transition from an event $ x $ to an event $ y $ when the former cannot influence the 
latter.

	Open and semiopen intervals of the real line are always denoted by square brackets, like 
$ ]0, 1[$, $[0, 1[, $ etc. The `extended numbers' $  \pm \infty $ are written $ \infty $ and $ \Oinf$, 
to avoid ambiguous expressions (see \cite{G4}). Marginal remarks and elementary proofs 
are written in small characters.

\subsection{Acknowledgements}\label{0.3}
I am indebted to my colleagues Ettore Carletti, Fernando Gliozzi and Nicola Pinamonti, 
for long discussions on the present subject.

\section{Enriched categories on a strict base}\label{s1}
	Categories enriched over monoidal categories are a classical topic of category theory 
\cite{EK, Ma, Ky}. The following simplified version is sufficient for our goals.

\subsection{Reviewing enriched categories}\label{1.1}
Informally, we recall that a category $ \C $ consists of {\em objects} $ X, Y, ... $ and 
{\em morphisms}, or {\em arrows}, $ f\c$ $X \to Y, ... $ between objects. The latter have a (partial) 
{\em composition law}: for each pair of consecutive morphisms $ f\c X \to Y $ and $ g\c Y \to Z $ 
there is a composed morphism $ gf\c X \to Z$. This partial operation is associative (whenever 
composition is possible) and every object $ X $ has an {\em identity}, written as $\id X\c X \to X$, 
or $ 1_X$, which acts as a unit for legitimate compositions.

\begin{small}

	We always assume that a category is {\em locally small}: the objects can form a class 
$\Ob \C$, but the morphisms $ X \to Y $ between two given objects always form a set 
$ \C(X, Y)$, called a {\em hom-set}. $ \C $ is {\em small} if the class $ \Ob \C $ is a set.

\end{small}

	We begin by recalling, or presenting, some elementary cases of `enriched categories'.

\Ndt (a) In a {\em preadditive category} $ \C$, each hom-set $ \C(X, Y) $ has a structure of 
abelian group, and the composition $ gf $ is additive in each variable
    \begin{equation}
g(f + f')  =  gf + gf',   \qq   (g + g')f  =  gf + g'f,
    \label{1.1.1} \end{equation}
(always assuming that these compositions are legitimate, of course). For instance, the category 
of real vector spaces (and linear mappings) has such a structure.

\Ndt (b) In an {\em ordered category} $ \C$, each hom-set $ \C(X, Y) $ has an order relation 
(reflexive, transitive and antisymmetric), and the composition is monotone in each variable
    \begin{equation}
f  \le f'  \implies   gf  \le gf',   \qq   g  \le g'   \implies   gf  \le g'f.
    \label{1.1.2} \end{equation}

	Let us note that here, because of transitivity, the previous condition is equivalent to 
global monotony: if $ f  \le f' $ and $ g  \le g'$, then $ gf  \le gf'  \le g'f'$.

\Ndt (c) {\em Terminology}. In case (a) we say that $ \C $ is a category {\em enriched over} 
$ \Ab$, the category of abelian groups (and homomorphisms), {\em equipped with the usual 
tensor product} $ A \te B$. We need the latter to express condition \eqref{1.1.1} by a 
homomorphism $ \C(X, Y) \te \C(Y, Z) \to \C(X, Z)$.

	In case (b) we say that $ \C $ is a category {\em enriched over} $ \Ord$, the category of 
ordered sets (and monotone mappings), equipped with the cartesian product $ A \ti B$. As 
we already remarked, condition \eqref{1.1.2} can be simply expressed by a monotone mapping 
$ \C(X, Y) \times \C(Y, Z) \to \C(X, Z) $ defined on the cartesian product.

\subsection{Categories}\label{1.2}
An ordinary (locally small) category can be viewed as enriched over $\Set$, the category of 
sets (and mappings), equipped with the cartesian product and its unit $E = \sing$, the 
singleton. We formalise this presentation, in view of a general definition of a category enriched 
over another.

	A ({\em locally small}) {\em category} $ \C $ consists of the following data:

\Ndt (a) a class $\Ob \C $ {\em of objects},

\Ndt (b) for each pair $ X, Y $ of objects, a set $ \C(X, Y) $ {\em of morphisms},

\Ndt (c) for each triple $ X, Y, Z $ of objects, a mapping {\em of composition}
$$k_{XYZ}\c \C(X, Y) \ti \C(Y, Z) \to \C(X, Z),$$

\Ndt (d) for each object $ X$, a mapping $ i_X\c \sing \to \C(X, X)$, {\em the identity} of $X$, 
giving the identity map $ \id X = i_X(*)$.

\skp These data have to satisfy three axioms, of associativity and unitarity: the following 
diagrams commute in $\Set$, for all the objects $ X, Y, Z, T $ of $ \C$ 
    \begin{equation} \begin{array}{c} 
    \xymatrix  @C=25pt @R=20pt
{
~\C(X, Y) \ti \C(Y, Z) \ti \C(Z, T)~   \ar[r]^-{1 \ti k}   \ar[d]_(.45){k \ti 1}    &
~\C(X, Y) \ti \C(Y, T)~   \ar[d]^(.45){k} 
\\ 
~\C(X, Z) \ti \C(Z, T)~   \ar[r]_-{k}  &   ~\C(X, T)~ 
}
    \label{1.2.1} \end{array} \end{equation}
    \begin{equation} \begin{array}{c} 
    \xymatrix  @C=25pt @R=18pt
{
~E \ti \C(X, Y)~   \ar[r]^-{\id \ti 1}   \ar[rd]_-{1}    &
~ \C(X, X) \ti \C(X, Y)~   \ar[d]^(.4){k}
\\ 
&   ~\C(X, Y)~ 
}
    \label{1.2.2} \end{array} \end{equation}
    \begin{equation} \begin{array}{c} 
    \xymatrix  @C=25pt @R=15pt
{
~\C(X, Y) \ti E~   \ar[r]^-{1 \ti \id}   \ar[rd]_-{1}    &
~\C(X, Y) \ti \C(Y, Y)~   \ar[d]^(.4){k} 
\\ 
&   ~\C(X, Y)~ 
}
    \label{1.2.3} \end{array} \end{equation}

	We are using the structural isomorphisms of the cartesian product in $ \Set$, that is the 
associativity isomorphism $ \al((a, b), c) = (a, (b, c)) $ and the unitarity isomorphisms 
$\la(*, a) = a = \rho(a, *)$.

\subsection{Monoidal categories}\label{1.3}
(a) Now we replace $ \Set $ with a {\em symmetric monoidal category} $ (\B, \sq, E)$, that 
will form the base of the enrichment.

	This is a (locally small) category $ \B $ equipped with a {\em binary product} $ A\sq B $ 
(often called the {\em tensor product}) that can be any functor $ \B \ti \B \to \B$, and a fixed 
object $ E $ (the {\em unit} of the product); furthermore, we have four natural isomorphisms 
of associativity, unitarity and commutativity
    \begin{equation} \begin{array}{c}
\al_{ABC}\c (A\sq B)\sq C \to A\sq (B\sq C),
\\[5pt]
\la_A\c E\sq A \to A,   \q   \rho_A\c A\sq E \to E,   \q   \si_{AB}\c A\sq B \to B\sq A.
    \label{1.3.1} \end{array} \end{equation}

	These isomorphisms are assumed to satisfy some axioms of coherence (cf.\ \cite{Ma}, 
Section VII.1), which will have no role here. A category $ \B $ with finite products has a 
{\em cartesian} structure $ (\B, \ti, \top) $ of this kind, produced by the cartesian product 
$A \ti B$, with unit given by the terminal object; the structural isomorphisms are a consequence 
of universal properties, and their coherence is automatically satisfied.

\Ndt (b) A symmetric monoidal category $ (\B, \sq , E) $ has an {\em associated forgetful functor}
    \begin{equation}
U\c \B \to \Set,   \qq   U(-)  =  \B(E, -),
    \label{1.3.2} \end{equation}
represented by the structural unit $ E$; it need not be faithful.

\Ndt (c) In the examples of Section \ref{1.1}, $ \Ab $ is a symmetric monoidal category, with 
respect to the usual tensor product $ A \te B$, with unit the additive group $ \bbZ$. $\Ord $ 
and $ \Set $ are cartesian monoidal categories. In all these cases, the associated forgetful 
functor is (isomorphic to) the usual one -- just the identity in the last case.

	The category $ \Cat $ of small categories (and functors) has also a cartesian monoidal 
structure; note that in this case the associated forgetful functor 
$ U = \Cat(\one, -)\c \Cat \to \Set $ (represented by the singleton category $ \one$) associates 
to a small category its set of objects, and is not faithful. (The reader may know that a 
category enriched on this base is a two-dimensional structure, called a 2-category.)

\subsection{Enrichment, in the strict case}\label{1.4}
(a) We assume now that the base $ (\B, \sq, E) $ is a {\em strict} symmetric monoidal category, 
the only case of interest here: this means that the product is strictly associative, unitary and 
commutative
    \begin{equation} \begin{array}{c}
A\sq (B\sq C) =  (A\sq B)\sq C,
\\[5pt]
E\sq A  =  A  =  A\sq E,   \q   A\sq B  =  B\sq A,
    \label{1.4.1} \end{array} \end{equation}
for all objects $A, B, C$. All the structural isomorphisms of the monoidal structure are 
thus identity morphisms of $ \B$, automatically coherent.

\Ndt (b) An {\em enriched category} $ \C $ over $ \B $ consists of the following data:

\Ndt - a class $\Ob \C $ {\em of objects},

\Ndt - for each pair $ X, Y $ of objects, a $\B$-object $ \C(X, Y) $ {\em of morphisms},

\Ndt - for each triple $X, Y, Z $ of objects, a $\B$--morphism {\em of composition}
$$k\c \C(X, Y)\sq \C(Y, Z) \to \C(X, Z),$$

\ndt - for each object $ X$, a $\B$--morphism $ i\c E \to \C(X, X)$, {\em the unit} of $X$.

\skp	These data are assumed to satisfy three axioms, of associativity and unitarity: the following 
diagrams have to commute in $ \B$, for all the objects $ X, Y, Z, T $ of $\C$ 
    \begin{equation} \begin{array}{c} 
    \xymatrix  @C=25pt @R=20pt
{
~\C(X, Y) \sq \C(Y, Z) \sq \C(Z, T)~   \ar[r]^-{1 \sq k}   \ar[d]_(.45){k \sq 1}    &
~\C(X, Y) \sq \C(Y, T)~   \ar[d]^(.45){k} 
\\ 
~\C(X, Z) \sq \C(Z, T)~   \ar[r]_-{k}  &   ~\C(X, T)~ 
}
    \label{1.4.2} \end{array} \end{equation}
    \begin{equation} \begin{array}{c} 
    \xymatrix  @C=30pt @R=15pt
{
~E \sq \C(X, Y)~   \ar[r]^-{\id \sq 1}   \ar[rd]_-{1}    &
~ \C(X, X) \sq \C(X, Y)~   \ar[d]^(.4){k} 
\\ 
&   ~\C(X, Y)~ 
}
    \label{1.4.3} \end{array} \end{equation}
    \begin{equation} \begin{array}{c} 
    \xymatrix  @C=30pt @R=15pt
{
~\C(X, Y) \sq E~   \ar[r]^-{1 \sq \id}   \ar[rd]_-{1}    &
~\C(X, Y) \sq \C(Y, Y)~   \ar[d]^(.4){k} 
\\ 
&   ~\C(X, Y)~ 
}
    \label{1.4.4} \end{array} \end{equation}

\begin{small}

\skp	The general case, where the base $ \B $ is only assumed to be a symmetric monoidal 
category, is obtained by inserting the structural isomorphisms \eqref{1.3.1}, in a rather 
obvious way.

\end{small}

\Ndt (c) Let us note that $ \C $ is not supposed to be a category. Rather, it has an 
{\em underlying category} $ \C_0 $ with the same objects and hom-sets 
$ \C_0(X, Y) = U(\C(X, Y))$, where $ U = \B(E, -)\c \B \to \Set $ is the canonical forgetful functor 
of the base.

\begin{small}

\skp	The composite of $ f\c E \to C(X, Y) $ and $ g\c E \to C(Y, Z) $ is given by
$$k(f\sq g)\c E = E\sq E \to C(X, Y)\sq C(Y, Z) \to C(X, Z),$$
and the identity of $ X $ is the unit $ i\c E \to \C(X, X)$. (We are implicitly using the fact that 
the functor $ U $ is lax monoidal, in a canonical way.) 

\end{small}

\Ndt (d) {\em If} $ U $ {\em is faithful}, we can think of $ \C $ as a category $ \C_0 $ where 
each hom-set $ \C_0(X, Y) $ is enriched with a structure $ \C(X, Y) $ in $ \B$, and we assume 
that the composition and identity mappings of $ \C_0 $ come out from the structure of $ \B$, as 
in (c). An example is presented below, in \ref{1.5}.

\begin{small}

\skp	This is also the case in several cases of enrichment over {\em weak} symmetric monoidal 
categories, as the examples of Section \ref{1.1}. A preadditive category can be simply presented 
as an ordinary category $ \C $ where each object $ \C(X, Y) $ is equipped with a structure of 
abelian group, and the (bilinear) composition mappings are homomorphisms 
$ \C(X, Y) \te \C(Y, Z) \to \C(X, Z)$. Similarly, an ordered category can be viewed as a 
category where each hom-set is an ordered set, and all composition mappings preserve the 
order relation.

\end{small}

\subsection{The truth-value base}\label{1.5}
As an elementary, useful example of a strict base, consider the ordinal category $ \two$, with 
two objects 0, 1, and a single non-identity arrow $ 0 \to 1$. It is equivalent to the category of all 
sets with at most one element, but it has a strict cartesian structure: 
$ 0 \ti 0 = 0 \ti 1 = 0$, $1 \ti 1 = 1$; the terminal is 1.

	A small category $ X $ enriched over $ \two $ amounts to a preordered set, that is, a set 
$X$, equipped with a truth-valued function $X(x, y) \in \{0, 1\}$, which we rewrite as $ x\prec y$ 
when $ X(x, y) = 1$; it has to satisfy the transitivity and reflexivity axioms
    \begin{equation}
(x \prec y, \;  y \prec z)   \implies   x \prec z,   \q   x \prec x   \qq   (\for  x, y, z \in X),
    \label{1.5.1} \end{equation}
which correspond to the composition arrow $ X(x, y) \ti X(y, z) \to X(x, z) $ and the identity arrow 
$ 1 \to X(x, x) $ of the enrichment.

	The (faithful!) forgetful functor $U = (1, -) $ embeds $\two$ as the full subcategory of 
$ \Set $ formed by $\es$ and a singleton. The category associated to $ X $ has the set 
$ X $ of objects, with one arrow $ x \to y $ precisely when $x \prec y$. This is, in fact, the 
usual way of presenting a preordered set as a (small) category.

\begin{small}

\skp	Even more basically, the singleton category $ \one $ has one object and its identity. 
A small category enriched on this base `is' just a set.

\end{small}

\section{Metric spaces as enriched categories}\label{s2}
	We now review Lawvere's generalised metric spaces, and their extension to real valued 
metrics, as enriched categories over a suitable strict base, studied in \cite{G4}. In both cases, 
the associated forgetful functor is far from being faithful, and of no interest.

\subsection{Lawvere's base}\label{2.1}
The strict monoidal category $ \Bw^+ = \OR^+ = [0, \infty] $ was introduced by Lawvere \cite{L1}, 
in 1974 (the original notation is $\bR$). Its objects are the extended positive numbers 
$ \la \in [0, \infty]$, with a unique morphism $ \la \to \mu $ when $ \la  \ge \mu $ (for the 
natural order), and a `tensor product' $ \la + \mu $ given by the extended sum, with unit $ 0$. 
All diagrams in $ \OR^+ $ commute (as in any category based on an ordered set).

	A small category $ X $ enriched over $ \OR^+ $ amounts to a generalised metric space, in 
the following sense: a set $ X $ equipped with hom-values $ X(x, y) \in [0, \infty]$, which we 
rewrite as $ \de(x, y)$, satisfying the triangle inequality and a `reflexivity axiom'
    \begin{equation}
\de(x, y) + \de(y, z)   \ge  \de(x, z),   \;\;\;\;   0   \ge  \de(x, x)   \q\;   (\for  x, y, z \in X).
    \label{2.1.1} \end{equation}

	Again, these laws represent the composition arrow $ X(x, y) + X(y, z)$ $\to X(x, z) $ and the 
identity arrow $ 0 \to X(x, x) $ of the enrichment. Note that the latter amounts here to 
$ \de(x, x) = 0$. The symmetry and separation axioms of classical metrics are not required.

\begin{small}

\skp	The forgetful functor $ U = \Bw^+(0, -)\c \Bw^+ \to \Set $ takes only two values: a singleton 
(at 0) and the empty set (elsewhere); the associated category $ X_0 $ is discrete, and the 
simplified view of enrichment does not work.

\end{small}

\subsection{Positive metrics and weighted algebraic topology}\label{2.2}
An influential part of the article \cite{L1}, developed by several authors, extends part of the 
theory of metric spaces to enriched categories, for instance Cauchy completeness.

	Here we are interested in another aspect: breaking the symmetry property, such a 
generalised metric space $ X $ can also express privileged directions. It was used in this sense 
by the present author, under the name of $\de$-{\em metric space}, or $\de$-{\em space}, as a 
basic structure for an enriched version of Directed Algebraic Topology \cite{G1, G3}, where 
the paths are measured and the forbidden ones are penalised with an infinite cost.

	The basic example is now the {\em standard} $\de$-{\em {\em line}} $ \de\bbR$, 
equipped with a $\de$-metric meant to hinder any backward move
    \begin{equation}
\de(x, x')  =  x' - x   \IF  x  \le x',   \q   \de(x, x')  =  \infty  \;  \text{ otherwise.}
    \label{2.2.1} \end{equation}

	Lipschitz paths and Lipschitz homotopies are parametrised on the {\em standard} 
$\de$-{\em interval} $ \de\bbI \sub \de\bbR$, equipped with the restricted $\de$-metric. 
The length of a Lipschitz path $ a\c \de\bbI \to X $ is given by the following least upper bound 
in $ \OR^+$:
    \begin{equation} \begin{array}{c}
L(a)  =  \sup_\Ut \, L_\Ut(a),
\\[5pt]
L_\Ut(a)  =  \Sum_j \, \de(a(t_{j-1}), a(t_j))    \q    (0 = t_0  \le t_1  \le ...  \le t_p =1),
    \label{2.2.2} \end{array} \end{equation}
where $\Ut = (t_j) $ stands for any finite increasing sequence in $ \bbI$.

	A $\de$-space $ X $ has a {\em fundamental weighted category} $\wPi_1(X) $ 
\cite{G2, G3}, an enrichment of the fundamental category of a directed space \cite{G1, G3}: it 
has for objects the points of $ X$, and for arrows $ \xi\c x \to x' $ the classes of Lipschitz paths 
$ a\c \de\bbI \to X $ from $ x $ to $ x'$, up to the equivalence relation generated by Lipschitz 
homotopies with fixed endpoints. The {\em weight} of such an arrow $ \xi\c x \to x' $ is defined 
as usual in a quotient, as the greatest lower bound of the length of the paths belonging to this 
equivalence class
    \begin{equation}
w(\xi)  =  \inf_{a\in\xi} \, L(a).
    \label{2.2.3} \end{equation}

	This gives $\wPi_1(X) $ a {\em weight}, or {\em generalised norm}, in $ \Bw^+ $ 
\cite{L1, G2, G3}. All this is extended to the more general framework of `spaces with weighted 
paths', whose finer quotients can interpret structures of Noncommutative Geometry: see 
\cite{G2}, Section 7, and \cite{G3}, Section 6.7.

\subsection{Real valued metrics}\label{2.3}
Following \cite{G4}, we now replace $ \Bw^+ $ by a larger structure $ \Bw $ based on the 
extended real line $ \OR= [\Oinf, \infty]$, to express `profits' and `losses', in any sense -- 
possibly related to a variable of any science, from Physics to Economy.

	This complete lattice is again viewed as a category, with morphisms $ \la  \ge \mu$, 
and endowed with a strict symmetric monoidal structure given by the extended sum 
$ \la + \mu$; here we let $ \Oinf + \infty = \infty$. (Adding any element should preserve the 
{\em initial} object $ \infty$, to get a `closed' monoidal structure: cf. \cite{G4}.)

	A {\em $\rho$-metric space}, or {\em $\rho$-space}, is a small category $ X $ enriched 
on $ \Bw$. This means a set $ X $ with a mapping $ \rho\c X \ti X \to \OR$ that satisfies the 
same conditions as $\de$ in \eqref{2.1.1}
    \begin{equation}
\rho(x, y) + \rho(y, z)   \ge  \rho(x, z),    \;\;\;\;   0 \ge  \rho(x, x)   \q\;   (\for x, y, z \in X),
    \label{2.3.1} \end{equation}
but the second law amounts here to $ \rho(x, x) = 0 $ {\em or} $ \Oinf$, for all $ x \in X $ 
(taking into account that $ \rho(x, x) + \rho(x, x)  \ge \rho(x, x))$.

	The extended real number $ \rho(x, y) $ expresses the cost of the transition from 
$ x $ to $ y$. Equivalently, the `antimetric' $ \ga(x, y) = - \rho(x, y) $ represents its gain 
(see \ref{3.1}). The transition is unaffordable, or forbidden, when its cost is $\infty$, and the 
gain is $\Oinf$.

	The `valuation', or `real-valued length', of a path $ a\c \bbI \to X $ (parametrised on 
the standard euclidean interval $ \bbI = [0, 1]$) is also an extended real number 
$ v_\rho(a)$, computed as in \eqref{2.2.2}
    \begin{equation}
v_\rho(a)  =  \sup_\Ut \, v_\Ut(a)  \in \OR,   \q   v_\Ut(a)  =  \Sum_j \, \rho(a(t_{j-1}), a(t_j)).
    \label{2.3.2} \end{equation}

	Two preprints of 1984 and 2015 also use metrics with values in $ \OR$, in 
other perspectives. Dealing with entropy, Lawvere introduced a real valued antimetric 
$ M(x, y) $ on the class of objects of a category $ \bbX $ endowed with an `entropy supply' 
$ s\c \bbX \to \OR$ (\cite{L2}, p.\ 4). Dealing with the Legendre-Fenchel transform for real 
vector spaces, Willerton used real valued metrics, called $\OR$-{\em metrics} \cite{Wi}.

\subsection{Complements}\label{2.4}
The basic example is now the {\em standard} $\rho$-{\em line} $ \rho\bbR$, with real metric
    \begin{equation}
\rho(x, y)  =  y - x.
    \label{2.4.1} \end{equation}

	This is actually a {\em linear} $\rho$-metric, which can be expressed as the variation 
$ \rho(x, y) = \Phi(y) - \Phi(x) $ of a function $ \Phi\c X \to \bbR$, the {\em potential} of $\rho$. 
Linear real metrics are important here (and trivial in the positive case); they could also be 
called `conservative metrics', being related to a potential.

	The asymmetry of a $\rho$-metric can be highlighted by two preorder relations, 
investigated in \cite{G4}, Section 4. Here we are mainly interested in the {\em reflective preorder}
    \begin{equation}
x \prec_\infty x'   \IF  \rho(x, x') < \infty,
    \label{2.4.2} \end{equation}
which means that the transition from $ x $ to $ x' $ is affordable. In $ \de\bbR $ and $ \de\bbI $ 
this gives the natural order.

\begin{small}

\skp	In $\rho\bbR$ this preorder is chaotic; on the other hand, the natural order corresponds 
to the {\em coreflective preorder} $x \prec_0 x'$, defined by $\rho(x', x)  \le 0$ \cite{G4}.

\end{small}

\section{Real antimetrics and Minkowski spacetime}\label{s3}
	Reversing inequalities, we briefly consider `antimetrics'. The theory is obviously equivalent 
to the previous one, but more natural in Relativity. The future timecone of a Lorentz vector 
space has a well-known `antinorm', which we use in \ref{3.4} to define an antimetric 
$ \ga(x, y) $ of the whole spacetime; the latter is $  \ge 0 $ if the event $ x $ can influence $ y$, 
and $ \Oinf $ otherwise. The geodesic paths have thus a {\em maximal} valuation with respect 
to $ \ga $ (the longest elapsed proper time) and {\em minimal} (if negative) with respect to the 
associated real metric $ \rho(x, y) = - \ga(x, y)$.

	We essentially follow the terminology of O'Neill's book on semi-Riemannian manifolds 
\cite{On}, although here it would be more natural to use the opposite signature of the metric tensor.

\subsection{Real antimetrics}\label{3.1}
(a) A (real valued) {\em antimetric} on a set $ X $ is a mapping $ \ga\c X \ti X \to \OR$ satisfying 
the axioms
    \begin{equation}
\ga(x, y) + \ga(y, z)   \le  \ga(x, z),    \;\;\;\;  0   \le  \ga(x, x)   \q\;   (\for x, y, z \in X),
    \label{3.1.1} \end{equation}

	The real value $ \ga(x, y) $ represents the {\em gain}, or {\em profit}, of the transition from 
$ x $ to $ y $ -- which is {\em unaffordable}, or {\em forbidden}, when $ \ga(x, y) = \Oinf$. The 
reversed triangular inequality expresses the fact that a combined transition should grant 
economies of scale; $ \ga(x, x) $ can only be $ 0 $ or $ \infty$.

\begin{small}

\skp	Formally, $ (X, \ga) $ is a small category enriched over the symmetric monoidal closed 
category $ \Bw^* $ associated to the complete lattice $ \OR$, equipped with its natural order 
and the extended sum where $ \infty + \Oinf = \Oinf $ (so that the initial object $ \Oinf $ is 
preserved by adding $\infty$).

\end{small}

\Ndt (b) The object $ (X, \ga) $ will also be called a $\ga$-{\em space}. The previous theory of 
$\rho$-spaces can be simply rewritten for these structures, reversing inequalities, 
$ \inf $ and $ \sup$, $\infty $ and $ \Oinf$.

	A {\em Lipschitz map} $ f\c X \to Y $ of $\ga$-spaces has a {\em Lipschitz constant} 
$ \la \in [0, \infty[ $ such that
    \begin{equation}
\la\ga_X(x, x') \le  \ga_Y(f(x), f(x'))   \q   (\for  x, x' \in X).
    \label{3.1.2} \end{equation}

	They form the category $ \ga_\infty\Mtr $ of $\ga$-{\em spaces and} {\em Lipschitz maps}. 
For $ \la = 1 $ this condition gives the wide subcategory $ \ga\Mtr $ of {\em gain-augmenting maps}.

	From \eqref{2.3.2}, a set-theoretical path $ a\c \bbI \to X $ has a {\em gain-valuation} 
$ v_\ga(a) $ defined as follows (with the usual meaning of a sequence $\Ut$, that is 
$0 = t_0  \le t_1  \le ...  \le t_p =1$)
    \begin{equation}
v_\Ut(a)  =  \Sum_j \, \ga(a(t_{j-1}), a(t_j)),   \q   v_\ga(a)  =  \inf_\Ut \, v_\Ut(a)  \in  \OR.
    \label{3.1.3} \end{equation}

\ndt (c) Let $ X $ be an abelian group, in additive notation. An antimetric $ \ga $ on $ X $ is 
{\em invariant up to translations} if and only if it can be expressed as $ \ga(x, y) = \Uga(y - x)$, 
by a (real) {\em antinorm} $ \Uga\c X \to \OR$, satisfying the axioms
    \begin{equation}
\Uga(x) + \Uga(y)   \le  \Uga(x + y),    \;\;\;\;   0   \le  \Uga(0)   \q\;   (\for x, y \in X).
    \label{3.1.4} \end{equation}

	Note that $ \Uga(0) $ can only be $ 0 $ or $ \infty$.

\subsection{The standard Lorentz vector space}\label{3.2}
(a) The $n$-dimensional {\em standard Lorentz vector space} $ V = \bbR^n_1 $ is a real vector 
space of dimension $ n  \ge 2$, equipped with the Lorentz scalar product
    \begin{equation}
\Bg(v, w)  =  \lan v, w \ran  =  \Sum_i \, \ep_i \, v_iw_i   \q   (i = 1, ..., n),
    \label{3.2.1} \end{equation}
of index  1 and signature $ \ep = (- 1, 1, ..., 1)$. $V $ can be viewed as the tangent vector space 
of the Minkowski spacetime, at any point.

	The scalar product is a non-degenerate symmetric bilinear mapping $ V \ti V \to \bbR$. 
We recall that `non-degenerate' means that $\lan v, w\ran$ only annihilates for all $ w $ if $ v = 0 $ 
(as shown letting $ w $ vary in the canonical base $ e_1, ..., e_n $ of $ \bbR^n$).

	The associated quadratic form
    \begin{equation}
\Bq(v)  =  \lan v, v\ran   =  \Sum_i \, \ep_i v_i^2,
    \label{3.2.2} \end{equation}
is indefinite; it will be used below to define an antinorm on $V$.

\Ndt (b) For a vector $ v \in V$, the number $ v_1 $ is the {\em time component} of $ v$, while 
the vector $\Vv = (v_2, ..., v_n) \in \bbR^{n-1}$ is its {\em spatial component}; we write as 
$ ||\Vv|| $ the euclidean norm of the latter, so that
    \begin{equation}
||\Vv||  =  (\Sum_{i>1} \, v_i^2)^{1/2},   \q  \Bq(v)  =  \lan v, v\ran   =  ||\Vv||^2 - v_1^2.
    \label{3.2.3} \end{equation}
%
%

\subsection{An antinorm for the Lorentz vector space}\label{3.3}
(a) The {\em closed future cone} of $ V $ is the following subset (related to the future causal 
cones of the spacetime, in \ref{3.4})
    \begin{equation}
\La  =  \{v \in V \sep \, \lan v, v\ran   \le 0,  \,  v_1  \ge 0\}  =  \{v \in V \sep  v_1  \ge ||\Vv||\},
    \label{3.3.1} \end{equation}
using the relation $ \lan v, v\ran  = ||\Vv||^2 - v_1^2$.

	Note that $ \la v + \mu w \in \La$, for all $ v, w \in \La $ and $ \la, \mu  \ge 0$. If 
$ v, w \in \La $ then $ \lan v, w\ran   \le 0$, as it follows applying the (ordinary) 
Cauchy-Schwarz inequality to the spatial part
$$\lan v, w\ran \, = \, (\Sum_{i>1} \, v_iw_i) - v_1w_1  \,  \le \,  ||\Vv||.||\Vw|| - v_1w_1 \,  \le \,  0.$$

\Ndt (b) The {\em open future cone} of $ V $ is its interior $ \La^\circ$, for the euclidean topology 
(and corresponds to the future timecones of \ref{3.4}.)
    \begin{equation}
\La^\circ  =  \{v \in V \sep \, \lan v, v\ran  < 0, \,  v_1 > 0\}  =  \{v \in X \sep  v_1 > ||\Vv||\}.
    \label{3.3.2} \end{equation}

	$\La^\circ$ is the {\em timecone} $ C(e_1)$, in the notation of \cite{On}, p.\ 143. The set 
of `timelike vectors', defined by the inequality $ \lan v, v\ran < 0$, is the disjoint union of 
$ \La^\circ$ with the opposite timecone $ - \La^\circ = C(- e_1)$.

\Ndt (c) The boundary of both cones is $ \dd\La = \La \setm \La^\circ$ (and corresponds to the 
future lightcones).

\Ndt (d) We define now a real antinorm on the vector space $V$, letting
    \begin{equation} \begin{array}{c}
\Uga(v)  = \sqrt{- \lan v, v\ran}  = \sqrt{v_1^2 -  ||\Vv||^2},
\q  \for \; v \in \La,
    \label{3.3.3} \end{array} \end{equation}
and $\Uga(v) =  \Oinf$ otherwise. All the vectors $v \notin \La $ are thus made unaffordable; 
for $ v \in \La$, the value $\Uga(v) \ge 0 $ is usually written as $|v|$, but we shall use 
this notation for all vectors.

\Ndt (e) We have to prove the reversed triangular inequality for $ \Uga$. On the cone 
$ \La^\circ$ this is an important, well-known fact (see \cite{On}, Corollary 5.31), related to 
the `reversed Schwarz inequality'. Its extension to $ V $ (once the value $ \Oinf $ is allowed) 
is easy.

\begin{small}

\skp	First, the inequality holds on the closed cone $ \La$: if $ v $ or $ w $ is in $ \dd\La$, 
taking into account that $ \lan v, w\ran  \le 0 $ and $ \Uga(v) \Uga(w) = 0$
    \begin{equation*} \begin{array}{c}
(\Uga(v + w))^2  \; = \; - \lan v + w, v + w\ran \;  = \;  - \lan v, v\ran  - 2\lan v, w\ran  - \lan w, w\ran
\\[5pt]
\q\;\;  \ge \;  \Uga(v)^2 + \Uga(w)^2 \; = \; (\Uga(v) + \Uga(w))^2.
    \label{} \end{array} \end{equation*}

	Finally, if one at least of $ v, w $ does not belong to $ \La$, we trivially have 
$ \Uga(v) + \Uga(w) = \Oinf  \le \Uga(v + w)$.

\end{small}

\subsection{An antimetric for the spacetime}\label{3.4}
(a) A {\em Minkowski spacetime} $ M $ is a time-oriented $n$-dimensional Lorentz manifold that 
is isometric to the Minkowski $n$-space $ X = \bbR^n_1$, for $ n  \ge 2 $ (slightly adapted from 
\cite{On}, Chapter 6, Definition 8).

	We shall identify $ M = X$, adopting an inertial coordinate system in $ M $ that preserves 
time-orientation (\cite{On}, Chapter 6, Definition 11). At any $ x \in X$, the tangent vector space 
$ T_xX $ can be identified with the Lorentz vector space $ V = \bbR^n_1 $ examined above, 
 (indefinite) scalar product $ \Bg(v, w) = \lan v, w\ran   = \Sum_i \, \ep_i \, v_iw_i$.

	$V$ and $X$ have thus the same underlying set, but a different structure:

\Ndt - $ V $ is a Lorentz vector space, and we still write its vectors as $ v, w, ...$,

\Ndt - $ X $ is a semi-Riemannian manifold, with Lorentz vector space $ T_xX = V $ 
at any point; its elements, written as $ x, y, ...$, are called {\em events}.

\skp	$V $ acts on the affine space $ X $ by the sum $ x + v $ (of $ V$). A difference $ y - x $ 
of events is interpreted as the vector $ v \in V = T_xX $ such that $ x + v = y$.

\Ndt(b) We recall that an event $ y \neq x $ can be reached by a material particle (resp.\  
a lightlike particle) starting at $ x $ if and only if it belongs to the {\em future timecone} 
$ x + \La^\circ$ (resp.\  the {\em future lightcone} $ x + \dd\La$). Globally, the event $ y $ 
can be influenced by $ x $ if and only if it belongs to the {\em future causal cone} 
$ C^+(x) = x + \La$. 

\Ndt (c) The antinorm $ \Uga\c V \to \OR$ defined in \ref{3.3}(d) gives an antimetric
    \begin{equation}
\ga\c X \ti X \to \OR,   \q   \ga(x, y)  =  \Uga(y - x).
    \label{3.4.1} \end{equation}

	Keeping $ x $ fixed, $ \ga(x, y) $ is finite and (weakly) positive when $ y $ belongs to the 
future causal cone $ C^+(x) = x + \La$, but takes value $ \Oinf $ on all events $ y $ which 
cannot be influenced by $ x$. It annihilates on the future lightcone $ x + \dd\La$.

\subsection{Lipschitz paths in the spacetime}\label{3.5}
Equivalently, the Minkowski spacetime $ M = X $ has a real valued metric, which takes values 
in $ ]\Oinf, 0] \cup \{\infty\}$
    \begin{equation}
\rho(x, y)  =  - \ga(x, y),
    \label{3.5.1} \end{equation}
although the gain function $ \ga$, which is {\em positive where it matters}, agrees better with 
the usual interpretation of spacetime.

	The reflective (pre)order in the $\rho$-space $ M $ is the {\em causality order}, defined by 
the future causal cones
    \begin{equation}
x  \le y  \; \iff \;   \rho(x, y) < \infty   \; \iff \;   y \in C^+(x) = x + \La.
    \label{3.5.2} \end{equation}

	It is easy to see that a mapping $ a\c \de\bbI \to M $ is a Lipschitz path if and only if it is 
weakly increasing for the reflective orders of these $\rho$-spaces
    \begin{equation}
t  \le t'  in  \bbI   \; \implies \;   a(t)  \le a(t')  \inn  M.
    \label{3.5.3} \end{equation}

	Special relativity is interested in the piecewise differentiable paths $ a\c \bbI \to M $ 
{\em that satisfy this condition}, and can thus represent the movement of a material or 
lightlike particle. Moreover, the integral of the metric tensor along the path measures the 
{\em proper duration} of the path itself -- formally zero for a light ray (which can bear no 
clock).
	
\subsection{Real valued metrics and topology}\label{3.6}
All this supports a claim of \cite{G4}: there is no sensible way of deriving a {\em definite} 
topology from a real-valued metric.

	The Minkowski spacetime $ M $ has a natural topology, defined by the manifold structure 
and homeomorphic to the euclidean space $ \bbR^n$. The `acceptable' paths $ a $ in $ M$, 
described above, are determined by the differentiable structure {\em and} the metric tensor (or
the associated $\rho$-metric); they are continuous for the natural topology, while all the 
topologies defined on $M$ by general procedures working on its $\rho$-metric (see \cite{G4}) 
are inadequate, or even trivial.

\section{Real antimetrics and Lorentz manifolds}\label{s4}
	More generally, any Lorentz manifold can be given a real valued antimetric.

\subsection{Lorentz vector spaces}\label{4.1} (a) Let $V $ be a real vector space of finite 
dimension.

	A {\em scalar product} on $ V $ is a non-degenerate symmetric bilinear mapping
    \begin{equation}
\Bg\c V \ti V \to \bbR,    \q    \Bg(v, v')  =  \lan v, v'\ran.
    \label{4.1.1} \end{equation}

	Then $ V $ admits an {\em orthonormal base} $ (e_i)_{i=1,..., n}$ such that
    \begin{equation}
\Bg(e_i, e_j)  =  \de_{ij} \ep_i,   \q   \ep_i  =  \Bg(e_i, e_i)  =   \pm 1,
    \label{4.1.2} \end{equation}
(with the Kronecker $ \de_{ij}$), as proved in \cite{On}, Lemma 24, p.\ 50. Then
    \begin{equation}
\Bg(u, v)  =  \Sum_i \, \ep_i \, u_iv_i   \q   (\for  u = \Sum_i \, u_i e_i, \;  v = \Sum_i \, v_i e_i).
    \label{4.1.3} \end{equation}

	The sequence $ (\ep_i)_{i=1,..., n}$ is called the {\em signature} of the scalar product 
with respect to our base $ (e_i)$, or more simply the {\em signature} of the base itself in $ V$. 
The number of elements $ e_i $ with $ \ep_i = - 1$ is the same for all the orthonormal bases 
(\cite{On}, Lemma 26, p.\ 51), and is called the {\em index} of $ \Bg$. If the index is 0, 
$ \Bg $ is positive definite, an {\em inner product}.

	Two vectors $ v, v' $ are said to be {\em orthogonal} if $ \lan v, v'\ran = 0$.

	The associated quadratic form is written as
    \begin{equation}
\Bq(v)  =  \lan v, v\ran.
    \label{4.1.4} \end{equation}

\Ndt (b) Let $ V $ be a {\em Lorentz vector space}, that is, a real vector space of dimension 
$ n  \ge 2$, equipped with a scalar product $ \Bg $ of index 1 (\cite{On}, p.\ 140). A 
non-trivial vector $ v \in V $ is said to be

\Ndt - {\em timelike} if $ \Bq(v) = \lan v, v\ran  < 0$,

\ndt - {\em lightlike} if $ \Bq(v) = \lan v, v\ran  =  0$,

\ndt - {\em spacelike} if $ \Bq(v) = \lan v, v\ran  >  0$.

\skp	The trivial vector is usually assumed to be spacelike, so that $ V $ is the disjoint union 
of the three subsets we are considering. If $ v $ is timelike, all its orthogonal vectors are 
spacelike (\cite{On}, Lemma 26, p.\ 141).

\subsection{Time-oriented Lorentz vector spaces}\label{4.2}
(a) By the previous remark, two timelike vectors $ u, u' $ have $ \lan u, u'\ran   \neq 0$. They are said to have the {\em same} (resp.\ {\em opposite}) {\em time-orientation} if $ \lan u, u'\ran  > 0 $ (resp.\ $ \lan u, u'\ran  < 0$).

	The Lorentz vector space $ V $ is {\em time-oriented} by the choice of a timelike vector $ u$. This determines the following `cones'
    \begin{equation} \begin{array}{lr}
\La  =  \{v \in V \sep \, \lan v, v\ran   \le 0, \,  \lan v, u\ran \ge 0\}  &  
\text{(the {\em causal cone} of  V)},
\\[5pt]
\La^\circ  =  \{v \in V \sep \, \lan v, v\ran  < 0, \,  \lan v, u\ran   \ge 0\}  &
\text{(the {\em timecone} of  V)},
\\[5pt]
\dd\La  =  \{v \in V \sep . \lan v, v\ran  = 0, \,  \lan v, u\ran   \ge 0\}  &
\text{(the {\em lightcone} of  V)}.
    \label{4.2.1} \end{array} \end{equation}

	It is important to note that any timelike vector $ u' $ with the same time-orientation gives 
the same cones, while any $ u' $ with opposite time-orientation gives the opposite cones 
$ - \La$, $ - \La^\circ$, $ - \dd\La $ (\cite{On}, p.\ 140, 146). The set of timelike vectors of $ V $ 
is the disjoint union $ \La^\circ \cup (- \La^\circ)$.

\Ndt (b) In particular $ \bbR^n_1$, the standard Lorentz vector space of dimension $ n$, is 
canonically time-oriented by the timelike vector $ e_1 = (1, 0, ..., 0) $ of the canonical base 
$ (e_i)$. Then any vector $ v = (v_i) $ has $ \lan v, e_1\ran  = v_1$, and the previous cones 
$ \La$, $ \La^\circ$, $ \dd\La $ coincide with the ones already introduced in \ref{3.3}.

\Ndt (c) The general case is essentially the same. In a time-oriented Lorentz vector space $ V$, 
let us fix an orthogonal base $ (e_i)_{i=1,..., n}$ such that $ \ep_1 = - 1 $ and the timelike vector 
$ e_1 $ agrees with the time-orientation of $ V $ (and determines it); the signature of the base 
is $ (\ep_i) = (- 1, 1, ..., 1)$.

	We have thus an isomorphism. defined on the standard Lorentz vector space of the same 
dimension
    \begin{equation}
f\c \bbR^n_1 \to V,
    \label{4.2.2} \end{equation}
that takes the canonical base of the latter to $(e_i)$. This isomorphism of vector 
spaces is an isometry (it preserves the scalar product, because of formula \eqref{4.1.3}), 
and preserves time-orientation. As a consequence, it preserves the quadratic form and the 
cones $\La, \La^\circ, \dd\La$.

\Ndt (d) By this isometry, the quadratic form $\Bq(v) = \lan v, v\ran$ gives a real antinorm on 
$ V $ (satisfying the reverse triangular inequality), as in \ref{3.3}
    \begin{equation} \begin{array}{lr}
|v|  =  (- \lan v, v\ran )^{1/2},  &  \for  v \in \La,
\\[5pt]
\;\;\;\;\;   =  \Oinf,   &   \text{otherwise.}
    \label{4.2.3} \end{array} \end{equation}
%
%

\subsection{Lorentz manifolds, path valuation and metric}\label{4.3}
(a) Let $ X $ be a {\em Lorentz n-manifold} of index 1, that is, a differentiable manifold of 
dimension $ n  \ge 2$, equipped with a tensor field $ \Bg$, which smoothly assigns to each 
point $ x \in X $ a scalar product $ \Bg_x = \lan -, -\ran _x $ of index 1, on the tangent vector 
space $ T_xX $ (\cite{On}, Definition 2, p.\ 55).

	 Moreover, we assume that $ X $ is {\em time-oriented}. For the sake of simplicity, this can 
be done by assigning a vector field $ U $ on $ X$, which is {\em timelike}, in the sense that 
every vector $ U_x $ is so: $ \Bg_x(U_x, U_x) < 0$, for all $ x \in M $ 
(see \cite{On}, Lemma 32, p.\ 145).

	We have thus, at each $ x \in M$, the causal cone of $ T_xX$
    \begin{equation}
\La_x =  \{v \in T_xX \sep \,  \lan v, v\ran _x  \le 0, \,  \lan v, U_x\ran _x  \ge 0\},
    \label{4.3.1} \end{equation}
and the real antinorm $ |v|_x $ defined in \eqref{4.2.3}.

	We say that a piecewise $C^1$ path is a {\em causal path} if, for each $ t \in \bbI$, the 
tangent vector $a'(t)$ belongs to the causal cone $\La_{a(t)}$. Of course, if $ a $ has an 
angular point at $ t \in \bbI$, we are assuming that both the left and right tangent vectors 
$a'_-(t) $ and $a'_+(t) $ belong to $ \La_{a(t)}$.

\Ndt (b) There is an additive valuation, restricted to the piecewise $ C^1 $ paths
    \begin{equation} \begin{array}{lr}
v_\Bg(a)  = \int_0^1 |a'(t)| \, dt,   \;\;  &  \text{for causal paths,}
\\[5pt]
\q\;\;   =  \Oinf,   &   \text{otherwise.}
    \label{4.3.2} \end{array} \end{equation}
which measures the proper elapsed time along the causal paths.

\Ndt (c) Plainly, $ v'_\Bg(a) = - v_\Bg(a) $ is also an additive valuation, restricted as above. 
We are interested in the $\rho$-metric associated to the latter
    \begin{equation}
\rho_\Bg(x, y)  =  \inf_a \, v'_\Bg(a)  =  - \sup_a \, v_\Bg(a),
    \label{4.3.3} \end{equation}
where $ a $ runs over the piecewise $ C^1 $ paths in $ X$, from $ x $ to $ y$. Thus 
$ \rho_\Bg(x, y) $ is a negative real number when $ x $ can influence $ y $ (i.e.\ there is a 
causal path from $ x $ to $ y$), and $ \infty $ otherwise; the antimetric 
$ \ga_\Bg(x, y) = - \rho_\Bg(x, y) $ measures the greatest proper elapsed time along such 
causal paths.

\Ndt (d) For the Minkowski spacetime $ X = \bbR^n_1$, considered in Section 2, one recovers 
the original $\rho$-metric
    \begin{equation}
\rho(x, y)  =  - \Uga(y - x),
    \label{4.3.4} \end{equation}
where $ y - x $ is a vector of the tangent vector space (at any point).

	On the other hand, the $\rho$-metric $ \rho' $ associated to $ v_\Bg $ (instead of 
$ v'_\Bg$) would give $ \rho'(x, y) = \inf_a \, v_\Bg(a) = 0 $ when $ x $ can influence $ y$, since 
there always is a piecewise $ C^1 $ light-line from $ x $ to $ y$, and the proper elapsed time 
along the latter is 0.



\begin{thebibliography}{99}


\bibitem[EK]{EK} S. Eilenberg and G.M. Kelly, Closed categories, in Proc. Conf. Categorical 
Algebra, La Jolla 1965, Springer, 1966, pp.\ 421--562. 

\bibitem[G1]{G1} M. Grandis, Directed homotopy theory, I. The fundamental category, 
Cah. Topol. G\'eom. Diff\'er. Cat\'eg. 44 (2003), 281--316. Available at

https://www.numdam.org/issues/CTGDC\_2003\_44\_4/

\bibitem[G2]{G2} M. Grandis, The fundamental weighted category of a weighted space 
(From directed to weighted algebraic topology), Homology Homotopy Appl. 9 (2007), 
221--256. Available at

https://intlpress.com/JDetail/1805807194991902721

\bibitem[G3]{G3}  M. Grandis, Directed Algebraic Topology, Models of non-reversible worlds, 
Cambridge Univ. Press, 2009. Available at

https://www.researchgate.net/publication/267089582

\bibitem[G4]{G4}  M. Grandis, Weighted algebraic topology, II (Metrics with real values), 
to appear.  

\bibitem[Ky]{Ky} G.M. Kelly, Basic concepts of enriched category theory, Cambridge 
University Press, 1982.

\bibitem[L1]{L1} F.W. Lawvere, Metric spaces, generalized logic and closed categories, 
Rend. Sem. Mat. Fis. Univ. Milano 43 (1974), 135--166. Republished in: Reprints 
Th. Appl. Categ. 1 (2002), 1--37. Available at

http://www.tac.mta.ca/tac/reprints/articles/1/tr1.pdf

\bibitem[L2]{L2} F.W. Lawvere, State Categories, Closed Categories, and the Existence 
Semi-Continuous Entropy Functions, IMA Preprint Series 84, University of Minnesota, 1984.

\bibitem[Ma]{Ma} S. Mac Lane, Categories for the working mathematician, Springer, 1971.

\bibitem[On]{On} B. O'Neill, Semi-Riemannian geometry, with applications to relativity, 
Academic Press, 1983.

\bibitem[Wi]{Wi} S. Willerton, The Legendre-Fenchel transform from a category theoretic 
perspective, arXiv 1501.03791v1, 2015.


\end{thebibliography}
\end{document}